\documentclass[a4paper, authoryear]{article}
\pdfoutput=1
\usepackage[parfill]{parskip}
\usepackage{natbib}
\usepackage[ruled,vlined,linesnumbered]{myalgorithm2e}
\SetCommentSty{textsf}
\usepackage{graphicx} 
\usepackage{placeins}
\usepackage{newproof}
\usepackage{setspace}
\usepackage{amsmath, amssymb}
\usepackage{threeparttable, longtable}
\usepackage[margin=1cm]{caption}
\usepackage[pdftitle={VPSolver 3: Multiple-choice Vector Packing Solver}, pdfauthor={Filipe Brandao}]{hyperref}
\usepackage{fullpage}
\usepackage{comment}

\newtheorem{corollary}{Corollary}

\newtheorem{property}{Property}
\newtheorem{definition}{Definition}

\newcommand{\vS}{\textsc{s}}
\newcommand{\vT}{\textsc{t}}
\newcommand{\vTs}{\textsc{t}^s}

\def\pauthor{
Filipe Brand\~ao\\
{\small INESC TEC and Faculdade de Ci\^{e}ncias,
Universidade do Porto, Portugal}\\
{\small\texttt{fdabrandao@dcc.fc.up.pt}}
}

\graphicspath{{figures_flow/}}

\begin{document}

\title{\textbf{VPSolver 3: Multiple-choice Vector Packing Solver}}
\author{\pauthor}
\date{\today}
\maketitle

\begin{abstract}
VPSolver is a vector packing solver based on an arc-flow formulation with graph compression.
In this paper, we present the algorithm introduced in VPSolver 3.0.0 for building compressed
arc-flow models for the multiple-choice vector packing problem.
\end{abstract}


\section{Introduction}

The vector bin packing problem (VBP), 
also called general assignment problem by
some authors, is a generalization of bin packing with multiple constraints. 
In this problem, we are
required to pack $n$ items of $m$ different types, represented by $p$-dimensional vectors, 
into as few bins as possible.
The multiple-choice vector bin packing problem (MVP) is a variant 
of VBP in which bins have several types (i.e., sizes and costs) 
and items have several incarnations (i.e., will take one of several possible sizes).

\cite{GeneralArcFlowPaper} present a general arc-flow formulation
with graph compression for vector packing.
This formulation is equivalent to the model of \cite{gomory1},
thus providing a very strong linear relaxation;
the largest absolute gap found in all the instances solved
in \cite{GeneralArcFlowPaper} was 1.0027.
Given a directed acyclic multi-graph containing every valid packing pattern
represented as a path from a source node to a target node,
the general arc-flow formulation is equivalent to the Gilmore-Gomory's model
with the same set of patterns as those represented as paths in the graph.

In \cite{GeneralArcFlowPaper}, a large variety of applications through reductions
to vector packing is presented.
\cite{MultipleChoiceVBP} extends this method to multiple-choice vector bin packing problems
by building a compressed arc-flow graph for each bin type.
In this method, a super source node is connected to the source node of each arc-flow graph and
each target node is connected to a super target.
Since every feasible packing pattern using any of the available bin types is represented
as a path from the super source to the super target, the general arc-flow formulation can be applied.

In this paper, we present an alternative algorithm for building compressed arc-flow graphs
for multiple-choice vector packing problems. This algorithm is a generalization of the 
direct Step-3 graph construction algorithm proposed in \cite{GeneralArcFlowPaper}.
If there is only one bin type and each item has a single incarnation, the new algorithm
produces exactly the same graph as the one produced by the original algorithm.
When multiple bin types exist, the algorithm produces on a single run
a graph containing all the valid packing patterns for all bin types.
This new algorithm is usually much faster than the method 
proposed in \cite{MultipleChoiceVBP}, and it usually produces smaller graphs.

The remainder of this paper is organized as follows.
Section~\ref{sec:arcflow} presents the arc-flow formulation for MVP.
In Section~\ref{sec:compression},
we show how to derive arc-flow models from dynamic programming recursions,
and how to obtain smaller models using graph compression.
Finally, Section~\ref{sec:algorithm} introduces the new graph construction and compression algorithm,
and Section~\ref{sec:conclusions} presents some conclusions.

\section{General arc-flow formulation for MVP}
\label{sec:arcflow}

For a given $i$, let $\mathbf{J}_i$ be the set of incarnations of item $i$, and
let $\mathbf{I} = \{(i,j) : i=1,\ldots,m,\ j \in \mathbf{J}_i\}$
be the set of items.
Let $(i,j)\in \mathbf{I}$ be the incarnation~$j$ of item~$i$,
and $w_{ij}$ its weight vector.
For the sake of simplicity, we define~$(0,0)$
as an item with weight zero in every dimension; this artificial item
is used to label loss arcs.
Let~$b_i$ be the demand of items of type~$i$, for $i=1,\ldots,m$.
Let~$b_0$ be the total number of items (i.e., $b_0 = \sum_{i=1}^{m} b_i$).
Let~$q$ be the number of bin types, and
let~$W_t$ and~$C_t$ be the capacity vector
and the cost of bins of type~$t$, respectively.

Given a directed acyclic multi-graph $G=(V,A)$ containing every valid packing 
pattern for each bin of type $t$ represented as a path from
the source $\vS$ to the target $\vT_t$,
and adding loss arcs connecting each target to the source,
the following arc-flow formulation can be used to model the corresponding multiple-choice
vector packing problem:
\begin{alignat}{3}
  & \mbox{minimize }   && \sum_{t=1}^q C_t f_{\vT_t\vS}^{0,0}\label{eq:new1}\\  
  & \mbox{subject to } \qquad && \sum_{(u,v,i,j)\in A:v=k} \hspace{-5mm}f_{uv}^{ij} \hspace{3mm}-\hspace{-3mm} \sum_{(u,v,i,j) \in A:u=k} \hspace{-5mm}f_{uv}^{ij} = 0\quad&& \mbox{for } k \in V,\label{eq:new2}\\    
      &&         & \sum_{(u, v, i, j) \in A:i=k} \hspace{-5mm}f_{uv}^{ij} \geq b_k, && k \in \{1,\ldots,m\} \setminus S, \label{eq:new3}\\ 
      &&         & \sum_{(u, v, i, j) \in A:i=k} \hspace{-5mm}f_{uv}^{ij} = b_k, && k \in S,  \label{eq:new4}\\
      &&         & f_{uv}^{ij} \leq b_i, && \forall (u,v,i,j) \in A, \label{eq:new5}\\
      &&         & f_{uv}^{ij} \geq 0, \mbox{ integer}, && \forall (u,v,i,j) \in A, \label{eq:new6}      
\end{alignat}
where each arc has four components~$(u, v, i, j)$
corresponding to an arc between nodes~$u$ and~$v$
associated to the incarnation $j$ of item $i$;
arcs~$(u,v,0,0)$ are loss arcs;
$f_{uv}^{ij}$~is the amount of flow along the arc~$(u, v, i, j)$;
$m$~is the number of different items;
$q$~is the number of bin types;
$b_i$~is the demand of items of type $i$;
and $S \subseteq \{1,\ldots,m\}$ is a subset of items whose demands 
are required to be satisfied exactly for efficiency purposes.
For having tighter constraints, one may set $S = \{i = 1,\ldots, m : b_i = 1\}$ 
(we have done this in our experiments). The only difference in relation
to the original arc-flow formulation for VBP is the fact that there is one target node for
each bin type instead of just one; on VBP instances this model is
exactly the same as the one proposed in \cite{GeneralArcFlowPaper}. 

In model (\ref{eq:new1})-(\ref{eq:new6}), the two types of demand constraints
and the upper bounds on variable values are mostly useful to take advantage of 
binary variables and multiple-choice constraints. 
A simplified version of this general arc-flow formulation that can also be used 
is the following:
\begin{alignat}{3}
  & \mbox{minimize }   && \sum_{t=1}^q C_t f_{\vT_t\vS}^{0,0}\label{eq:simple1}\\  
  & \mbox{subject to } \qquad && \sum_{(u,v,i,j)\in A:v=k} \hspace{-5mm}f_{uv}^{ij} \hspace{3mm}-\hspace{-3mm} \sum_{(u,v,i,j) \in A:u=k} \hspace{-5mm}f_{uv}^{ij} = 0,\quad&& \mbox{for } k \in V,\label{eq:simple2}\\    
      &&         & \sum_{(u, v, i, j) \in A:i=k} \hspace{-5mm}f_{uv}^{ij} \geq b_k, && k = 1,\ldots,m, \label{eq:simple3}\\ \
      &&         & f_{uv}^{ij} \geq 0, \mbox{ integer}, && \forall (u,v,i,j) \in A, \label{eq:simple6}      
\end{alignat}

The following pattern based formulation generalizes \cite{gomory1}
formulation to the multiple-choice vector packing problem and it is equivalent
to model (\ref{eq:new1})-(\ref{eq:new6}).
For bins of type $t$,
let column vectors 
$a_{kt} = \langle a_{kt}^{ij} \rangle_{(i,j) \in \mathbf{I}}$
represent all feasible packing patterns~$k$; 
each element~$a_{kt}^{ij}$ represents the number of times incarnation $j$ of item $i$ is used in the pattern.
Let~$x_{kt}$ be a decision variable that
designates the number of times pattern~$k$ using bins of type $t$ is used.  
The MVP can be modeled in terms of these variables as follows:
\begin{alignat}{3}
  & \mbox{minimize }          && \sum_{t=1}^{q} \sum_{k \in K_t} C_t x_{kt} \label{eq:pat1}\\  
  & \mbox{subject to } \qquad && \sum_{t=1}^{q} \sum_{k \in K_t} \sum_{j \in  \mathbf{J}_i} a_{kt}^{ij} x_{kt} \geq b_i, \qquad && i = 1,\ldots,m, \label{eq:pat2}\\
  &                           && x_{kt} \geq 0, \mbox{ integer},    \qquad && t=1,\ldots,q,\ \forall k \in K_t, \label{eq:pat3}
\end{alignat}
where every valid packing pattern $k \in K_t$ for bins of type $t$ satisfies:
\begin{alignat}{3}
    & \sum_{(i,j) \in  \mathbf{I}} a_{kt}^{ij}  w_{ij}^d \leq W^{d}_t,\qquad &&d=1,\ldots,p, \\
    & a_{kt}^{ij} \geq 0, \mbox{ integer}, && (i,j) \in \mathbf{I}.
\end{alignat}

Since there is flow conservation in every node
in models (\ref{eq:new1})-(\ref{eq:new6}) and (\ref{eq:simple1})-(\ref{eq:simple6}),
the arc-flow solutions are circulations.

\begin{corollary}[Flow decomposition for circulations]\label{flowdecomposition}
Any non-negative feasible circulation flow can be decomposed
into the sum of flows around directed cycles.
\end{corollary}

Corollary~\ref{flowdecomposition} follows directly from the Flow Decomposition Theorem (see, e.g., \citealt{ahuja-magnanti-orlin-93}).

\begin{property}[Equivalence to the pattern-based model]\label{equivalence_pattern}
For a graph with all valid packing patterns represented as paths 
from $\vS$ to $\vT_t$ for each bin type $t$, 
model (\ref{eq:simple1})-(\ref{eq:simple6})  is equivalent to the pattern-based
model (\ref{eq:pat1})-(\ref{eq:pat3}) with the same set of patterns.
\end{property}
\begin{proof}  
We apply Dantzig-Wolfe decomposition to model~(\ref{eq:simple1})-(\ref{eq:simple6})
keeping (\ref{eq:simple1}) and (\ref{eq:simple3}) in the master problem, and
(\ref{eq:simple2}) and (\ref{eq:simple6}) in the subproblem. As the subproblem
is a totally unimodular flow model whose solutions can be decomposed into cycles
(each including one of feedback arcs associated with a bin type),
only valid packing patterns are generated, and
we can substitute (\ref{eq:simple2}) and (\ref{eq:simple6}) 
by the patterns and obtain a pattern-based model.  
From this equivalence follows that lower bounds
provided by both models are the same when the same set of patterns is considered,
and the solution space is exactly the same.
\end{proof}

\begin{property}[Equivalence to the pattern-based model]\label{equivalence_pattern}
For a graph with all valid packing patterns represented as paths 
from $\vS$ to $\vT_t$ for each bin type $t$, 
model (\ref{eq:new1})-(\ref{eq:new6})  is equivalent, in terms of lower bound at the root node, 
to the pattern-based model (\ref{eq:pat1})-(\ref{eq:pat3}) with the same set of patterns.
\end{property}
\begin{proof}  
We apply Dantzig-Wolfe decomposition to model~(\ref{eq:new1})-(\ref{eq:new6})
keeping (\ref{eq:new1}), (\ref{eq:new3}) and (\ref{eq:new4}) in the master problem and
(\ref{eq:new2}), (\ref{eq:new5}) and (\ref{eq:new6}) in the subproblem. As the subproblem
is a totally unimodular flow model whose solutions can be decomposed into cycles
(each including one of feedback arcs associated with a bin type),
only valid packing patterns are generated, and
we can substitute (\ref{eq:new2}), (\ref{eq:new5}) and (\ref{eq:new6}) 
by the patterns and obtain a pattern-based model.  
From this equivalence follows that lower bounds
provided by both models are the same when the same set of patterns is considered.
The equality constraints (\ref{eq:new4}) and the upper bound on 
variable values (\ref{eq:new5})
have no effect on the lower bounds, 
since we are only excluding solutions that satisfy the demand of some items with excess, 
and for these solutions there are equivalent solutions
that do not exceed the demand (recall that every valid packing pattern is represented in the graph).
\end{proof}

\section{Arc-flow models and graph compression}
\label{sec:compression}

For cutting and packing problems,
arc-flow models equivalent to pattern-based models can be easily derived 
from the dynamic programming recursion of the underlying knapsack subproblem 
(see, e.g., \citealt{Wolsey1977527}).
The main challenge is to find a compact representation of the patterns
in a reasonable amount of time.
In Section~\ref{sec:deriving}, we show how the arc-flow models
can be easily derived from dynamic programming recursions, and in Section~\ref{sec:graphcomp},
we show how to obtain smaller arc-flow models using graph compression.

\subsection{Deriving arc-flow models from dynamic programming recursions\label{sec:deriving}}

\begin{definition}[Order]\label{def:order}
Items are sorted in decreasing order by the sum of normalized weights
($\alpha_{ij} = \sum_{d=1}^{p} w_{ij}^d/\max\{W^d_t : t=1,\ldots,q\}$), using decreasing
lexicographical order in case of a tie.
\end{definition}

The multiple-choice vector packing problem has an underlying knapsack subproblem
on the capacity constraints of each bin.
Let $\mathbf{I}'$ be the set of items sorted according Definition~\ref{def:order}.
Let $\pi_i$ be the profit of item $i$.
For the sake of simplicity, let $w_{k} = w_{ij}$, $b_k = b_i$, $\pi_k = \pi_i$, and arcs $(u, v, k) = (u, v, i, j)$, for $k=1,\ldots,|\mathbf{I}'|$ and $(i,j) = \mathbf{I}'_k$.
A dynamic programming recursion for the knapsack subproblem can be defined as follows:

\begin{equation*}
\begin{array}{rll}
D(x, k, c) &=& \left\{\begin{array}{ll}
\min\{C_t : t=1,\ldots,q,\ x \leq W_t \} & \mbox{if } k=|\mathbf{I}'|+1,\\
\min\{D(x+w_{k}, k, c+1)-\pi_{k}, D(x, k+1, 0) \} & \mbox{if }  c+1 \leq b_k \mbox{ and } \exists t\ x+w_{k} \leq W_t,\\
D(x, k+1, 0) & \mbox{otherwise}.\\
\end{array}
\right.\\
\end{array}
\end{equation*}

Note that this dynamic programming recursion could be used to find the most attractive column
in a column generation algorithm based on model~(\ref{eq:pat1})-(\ref{eq:pat3}).
At each iteration of the column-generation process, 
a subproblem is solved and a column (pattern) is introduced
in the restricted master problem if its reduced cost is strictly less than zero.
Let $\pi_k$ be the shadow price of the demand constraint associated with item $k$.
The reduced cost of the most attractive pattern is given by $D(\langle 0 \rangle_{d=1,\ldots,p}, 1, 0)$.
In our method, instead of using column-generation in an iterative process, we construct a
graph containing every valid packing pattern, and this graph can be derived from 
this dynamic programming recursion.

From this dynamic programming recursion $D$, we can easily derive an arc-flow model.
Consider each dynamic programming state as a node, and each recursive call as an arc.
The graph can be obtained as follows.
Let $G(x, k, c) = \{((x, k, c), \vTs_t, 0): t = 1,\ldots,q \} \mbox{ if } k=|\mathbf{I}'|+1$,
$G(x, k, c) = \{((x, k, c), (x+w_{k}, k, c+1), k), ((x, k, c), (x, k+1, 0), 0)\} \cup G(x+w_{k}, k, c+1) \cup G(x, k+1, 0) \mbox{ if }  c+1 \leq b_k \mbox{ and } \exists t\ x+w_{k} \leq W_t$,
and $G(x, k, c) = G(x, k+1, 0)$ otherwise. The source node~$\vS$ is $(\langle 0 \rangle_{d=1,\ldots,p}, 1, 0)$, 
the set of target nodes is $T=\{ \vTs_t :t=1,\ldots,q \}$ (the base cases of the dynamic programming recursion),
the set of arcs is $A = G(\langle 0 \rangle_{d=1,\ldots,p}, 1, 0)$,
and the set of vertices is $V = \{u : (u, v, k) \in A \} \cup \{v: (u, v, k) \in A \}$.

The dynamic programming recursion $D$ is equivalent to the following totally unimodular flow
problem:
\begin{alignat}{3}
  & \mbox{minimize }   && \sum_{t=1}^{q} \sum_{(u,v,k) \in A: v=\vTs_t} \hspace{-6mm}C_t f_{uv}^k\hspace{3mm}  - \sum_{(u,v,k) \in A: k \neq 0} \hspace{-5mm}\pi_k f_{uv}^{k}\hspace{3mm}\\  
  & \mbox{subject to } \qquad && \sum_{(u,v,k)\in A:u=v'} \hspace{-5mm}f_{uv}^{k} \hspace{3mm}-\hspace{-3mm} \sum_{(u,v,k) \in A:v=v'} \hspace{-5mm}f_{uv}^{k} = 0, \quad && \mbox{for } v' \in V \setminus (\{\vS\} \cup T),\\   
  & && \sum_{(u,v,k) \in A:u=\vS} \hspace{-5mm}f_{uv}^{k}\hspace{3mm} = 1, \quad &&\\       
      &&         & f_{uv}^{k} \geq 0, && \forall (u,v,k) \in A.
\end{alignat}

The dual of this flow problem,
which corresponds exactly to the dynamic programming recursion defined above,
is the following:
\begin{alignat}{3}
  & \mbox{maximize }   && \theta'_{\vS}\\  
  & \mbox{subject to } \qquad && \theta'_{u} \leq C_t, \quad && \mbox{for } (u,v,k) \in A, k = 0, \exists t\ v=\vTs_t,\\
  & && \theta'_{u} \leq \theta'_{v}, \quad && \mbox{for }  (u,v,k) \in A, k = 0, v \notin T,\\
  & && \theta'_{u} \leq \theta'_{v}-\pi_{k}, \quad && \mbox{for }  (u,v,k) \in A, k \neq 0.       
\end{alignat}
This relationship between dynamic programming recursions for knapsack problems and arc-flow models
is one of the results of \citealt{Wolsey1977527},
which leads to a natural method for obtaining arc-flow formulations for cutting 
and packing problems.

\subsection{Graph compression\label{sec:graphcomp}}

In the graphs derived from the dynamic programming recursion $D$,
each internal node is connected to at most two other nodes:
a node in its level (using the current item), and another in the level above (not using the current item).
This graph can be seen as the Step-2 graph of \cite{GeneralArcFlowPaper}
as it already divides the graph into levels, one for each item incarnation, and therefore breaks symmetry.

By adding the feedback loss arcs connecting each target node to the source (i.e., $\{(\vTs_t, \vS, 0): t=1,\ldots,q\}$),
the arc-flow model derived directly from the dynamic programming recursion can be used with  models (\ref{eq:new1})-(\ref{eq:new6}) or (\ref{eq:simple1})-(\ref{eq:simple6})
to solve the corresponding multiple-choice vector packing problem. However,
since there is a constraint for each node and a variable for each arc, the model size may be a problem,
and therefore this can only be done for very small instances. 
Graph compressing is used in order to solve this problem as it usually allows us to obtain
reasonably small graphs even when the straightforward approach would lead 
to models with many millions of variables and constraints.

Graph compression consists of relabelling the graph.
In each compression step, a new graph $G'~=~(V',A')$ is constructed
by creating a set of vertices $V'=\{\phi(v) : v \in V \}$ and
a set of arcs $A' = \{(\phi(u),\phi(v),k) : (u,v,k) \in A, \phi(u) \neq \phi(v) \}$,
where~$\phi$ is the map between the original and new labels.
This relabelling procedure must assure that no valid packing pattern is removed,
and that no invalid packing pattern is added.

The main compression step is applied to the \mbox{Step-2} graph (i.e., a graph with level dimension labels which break symmetry).
In the \mbox{Step-3} graph, the longest paths to the target in each dimension are used
to relabel the internal nodes ($V \setminus (\{\vS\} \cup T)$), dropping the level dimension labels (i.e., $k$ and $c$) of each node.
Let $\langle \varphi^d(u) \rangle_{d=1,\ldots,p}$ be the new label of 
node~$u$ in the \mbox{Step-3} graph, where
\begin{alignat}{3}
\varphi^d(u) & = & \left\{ \begin{array}{ll}                
                W^d_t & \mbox{if }  u = \vTs_t \mbox{ (base case)}, \\
                \min_{(u',v,k) \in A:u'=u}\{\varphi^d(v) - w_k^d\}  & \mbox{otherwise.}\\
                \end{array}\right.           
\end{alignat}
In the paths from $\vS$ to $\vTs_t$ in the \mbox{Step-2} graph usually there is slack in some dimension.
In this process, we are moving this slack as much as possible to the beginning of the paths.

Finally, in the final compression step, a new graph is constructed once more. 
In order to try to reduce the graph size even more, 
we relabel the internal nodes once more using the longest paths from the source in each dimension.
Let $\langle \psi^d(v) \rangle_{d=1,\ldots,p}$ be the label of node~$v$ in the
\mbox{Step-4} graph, where
\begin{alignat}{3}
\psi^d(v) & = & \left\{ \begin{array}{ll}
                0 & \mbox{if }  v = \vS \mbox{ (base case)}, \\
                \max_{(u,v',k) \in A:v'=v}\{\psi^{d}(u) + w_k^d\}  & \mbox{otherwise.}\\
                \end{array}\right.           
\end{alignat}

\section{Graph construction and compression algorithm}
\label{sec:algorithm}

Algorithm~\ref{alg:build} is a generalization of the algorithm proposed in \cite{GeneralArcFlowPaper}.
It builds the \mbox{Step-3} graph directly
in order to avoid the construction of huge \mbox{Step-1} and \mbox{Step-2} graphs
that may have many millions of vertices and arcs.
We start by sorting the items according the order defined in Defenition~\ref{def:order} in line~\ref{build:sort}.
Algorithm~\ref{alg:build} uses dynamic programming to build the \mbox{Step-3} graph
recursively over the structure of the \mbox{Step-2} graph
(without building it).
The basic idea for this algorithm comes from the fact that, in the main compression step, the label of any internal node ($\varphi(u) = \langle \min\{\varphi^d(v) - w_{ij}^d : (u',v,i,j) \in A, u'=u\} \rangle_{d=1,\ldots,p}$) only depends
on the labels of the two nodes to which it is connected; 
a node in its level (line~\ref{alg:case2}) and another in the level above (line~\ref{alg:case1}).
After directly building the \mbox{Step-3} graph from the instance's 
data using this algorithm, we 
apply the final compression step (line~\ref{build:finalcomp}) using Algorithm~\ref{alg:finalcomp},
and connect the internal nodes to the targets (line~\ref{build:contgts}) using Algorithm~\ref{alg:contgts}.
Since parallel arcs associated to the same item type but different incarnations are redundant,
all redundant arcs are removed in line~\ref{build:parallel}.
In practice, this method allows obtaining arc-flow models even for large benchmark instances quickly.

\begin{algorithm}[!h]
\caption{Graph construction and compression algorithm}\label{alg:build}
\SetEndCharOfAlgoLine{}
\SetKwInOut{Input}{input}
\SetKwInOut{Output}{output}

\SetKwFunction{Build}{build}
\SetKwFunction{BuildGraph}{buildGraph}
\SetKwFunction{ConnectTargets}{connectTargets}
\SetKwFunction{FinalCompression}{finalCompression}
\SetKwFunction{Sort}{sorted}
\SetKwFunction{Key}{key}
\SetKwFunction{Len}{len}
\SetKwFunction{Reverse}{reverse}
\SetKwFunction{Lift}{lift}
\SetKwFunction{HPos}{highestPosition}
\SetKwData{varitems}{items}
\SetKwData{vardp}{dp}
\SetKwData{varx}{x}
\SetKwData{vari}{i}
\SetKwData{varmx}{mx}
\SetKwData{varmi}{mi}
\SetKwData{NIL}{NIL}

\SetInd{0.5em}{1.5em}

\SetKwBlock{Function}{}{}

\Input{$\mathbf{I}$ - set of items; $w$ - item weights; $b$ - demands; $q$ - number of bin types; $W$ - capacity vectors}
\Output{$V$ - set of vertices; $A$ - set of arcs; $\vS$ - source; $\vTs$ - targets}

\textbf{function} $\BuildGraph(\mathbf{I}, w, b, q, W)$:
\Function{
$\vardp[x, k, c] \gets \NIL, \mbox{ \bf for all } x, k, c$\tcp*[r]{dynamic programming table}
$\mathbf{I}' \gets \Reverse(\Sort(\mathbf{I}, \Key=\lambda (i, j) . (\sum_{d=1}^{p} w_{ij}^d/\max\{W_t^d : t=1,\ldots,q \}, w_{ij})))$\tcp*[r]{sort item incarnations\label{build:sort}}
\textbf{function} $\Lift(x, k, c)$: \tcp*[f]{auxiliary function: lift \vardp states solving knapsack/longest-path problems in each dimension}
\Function{
	\Input{$x$ - used capacity; $k$ - current item; $c$ - number of times the current has been used}
    \textbf{function} $\HPos(d, t)$:\label{alg:hpos}
    \Function{
	   \Return $\min \left \{  
	    W_t^d-\sum_{j = k}^{|\mathbf{I}'|} w_{\mathbf{I}'_{j}}^d y_j :
	     \begin{array}{ll}
	     \sum_{j = k}^{|\mathbf{I}'|} w_{\mathbf{I}'_{j}}^d y_j \leq W_t^d-x^d,\\
	     y_k \leq b_k-c,\ y_j \leq b_j,\ j=k+1,...,|\mathbf{I}'|,\\
	     y_j \geq 0, \mbox{ integer},\  j=k,...,|\mathbf{I}'|  
	     \end{array}
		\right \}       	   
	   $
    }
    \Return $\left \langle \min\{\HPos(d, t): t =1,\ldots,q, x \leq W_t\} \right \rangle_{d=1,\ldots,p}$    
}
$(V, A) \gets (\{\ \}, \{\ \})$\\
\textbf{function} $\varphi(x, k, c)$:
\Function{
\Input{$x$ - used capacity; $k$ - current item; $c$ - number of times the current has been used}
$x \gets \Lift(x,k,c)$\tcp*[r]{lift $x$ in order to reduce the number of \vardp states}\label{alg:lift}
\If(\tcp*[f]{avoid repeating work }){$\vardp[x,k,c] \neq \NIL$}{
    \Return $\vardp[x,k,c]$
}
$u \gets \left \langle \min\{W_t^d : t=1,\ldots,q, x \leq W_t \} \right \rangle_{d=1,\ldots,p}$\tcp*[r]{base case of $\varphi(x,k,c)$ if there are no arcs leaving the node}
\If(\tcp*[f]{option 1: do not use the current item (go to the level above)}\label{alg:case1}){$k < |\mathbf{I}'|$}{
    $up_x \gets \varphi(x, k+1, 0)$\\
    $u \gets up_x$\tcp*[r]{value of $\varphi(x,k,c)$ if no more items of the current type are introduced}
}
$(i, j) \gets \mathbf{I}'_k$\\
\If(\tcp*[f]{option 2: use the current item}\label{alg:case2}){$c < b_i \mbox{ \bf and }  x+w_{ij} \leq W_t \mbox{ \bf for any } t=1,\ldots,q$}{    
    $v \gets \varphi(x+w_{ij}, k, c+1)$\\
    $u \gets \left \langle \min(u^d, v^d-w_{ij}^d) \right \rangle_{d=1,\ldots,p}$\tcp*[r]{update the value of $\varphi(x,k,c)$}
    $A \gets A \cup \{(u, v, i, j)\}$\tcp*[r]{connect $u$ to the node resulting from option 2}
    $V \gets V \cup \{u, v\}$
}
\If(){$k < |\mathbf{I}'|${ \bf and }$u \neq up_x$}{
    $A \gets A \cup \{(u,up_x,0,0)\}$\tcp*[r]{connect $u$ to the node resulting from option 1}
    $V \gets V \cup \{up_x\}$
}
$\vardp[x,k,c] \gets u$

\Return $u$\tcp*[r]{returns $u=\varphi(x,k,c)$}
}
$\vS \gets \varphi(x=\left \langle 0 \right \rangle_{d=1,\ldots,p},k=1,c=0)$\tcp*{build the graph}
$(V, A, \vS) \gets \FinalCompression(V, A, \vS, w)$\tcp*{final compression step\label{build:finalcomp}}
$(V, A, \vS, \vTs) \gets \ConnectTargets(V, A, \vS, q, W)$\tcp*{connect the internal nodes to the targets\label{build:contgts}}
$A \gets \{(u, v, i, j) \in A : (u, v, i, j') \notin A, \forall j' < j\}$\tcp*{remove parallel arcs associated to the same item type\label{build:parallel}}
\Return $(G = (V, A), \vS, \vTs)$
}
\end{algorithm}

The dynamic programming states are identified by
the space used in each dimension~($\langle x^d \rangle_{d=1,\ldots,p}$), the current item~($k$)
and the number of times it has already been used~($c$).
In order to reduce the number of states, we lift (line~\ref{alg:lift})
each state by solving (using again dynamic programming though this is not explicit in the algorithm) knapsack/longest-path problems 
in each dimension considering the remaining items (line~\ref{alg:hpos});
we try to increase the space used in each dimension
to its highest value considering the valid packing patterns for the remaining items.
Note that all valid bin sizes must be considered in the lifting procedure.

If there is just one bin type, this algorithm works exactly as the original one.
When there are multiple bin types, the main difference are the lift procedure,
which needs to take in consideration all valid bin sizes, 
and the connection of internal nodes to the targets.
When there are multiple valid targets for the same node, we need to connect the node
to each of them, or take advantage of a transitive reduction to connect each node to
as little targets as possible (as we do in Algorithm~\ref{alg:contgts}). For instance, in the variable bin size problem,
since there is just one dimension, the transitive reduction allows us to connect each internal node
to just one target (i.e., it uses only 1 additional arc per node instead of $q$ additional arcs).

\begin{algorithm}[!h]
\caption{Final compression step}\label{alg:finalcomp}
\SetEndCharOfAlgoLine{}
\SetKwInOut{Input}{input}
\SetKwInOut{Output}{output}

\SetKwFunction{Build}{build}
\SetKwFunction{BuildGraph}{buildGraph}
\SetKwFunction{FinalCompression}{finalCompression}
\SetKwFunction{Sort}{sorted}
\SetKwFunction{Key}{key}
\SetKwFunction{Len}{len}
\SetKwFunction{Reverse}{reverse}
\SetKwFunction{Lift}{lift}
\SetKwFunction{HPos}{highestPosition}
\SetKwData{varitems}{items}
\SetKwData{vardp}{dp}
\SetKwData{varx}{x}
\SetKwData{vari}{i}
\SetKwData{varmx}{mx}
\SetKwData{varmi}{mi}
\SetKwData{NIL}{NIL}

\SetInd{0.5em}{1.5em}

\SetKwBlock{Function}{}{}

\Input{$V$ - set of vertices; $A$ - set of arcs; $\vS$ - source; $w$ - item weights}
\Output{$V$ - set of vertices; $A$ - set of arcs; $\vS$ - source}

\textbf{function} $\FinalCompression(V, A, \vS, w)$:
\Function{
$\psi(\vS) \gets \left \langle 0 \right \rangle_{d=1,\ldots,p}$\\
\ForEach(\tcp*[f]{for each vertex in reverse topological order of the transpose graph}){$v \in\Sort(V \setminus \{\vS\})$}{
	$\psi(v) \gets \left \langle \max \{\psi^d(u) + w_{ij}^d\ : (u, v', i, j) \in A,\ v'=v\} \right \rangle_{d=1,\ldots,p}$\\
}
$\vS \gets \psi(\vS)$\\
$V \gets \{\psi(u) : u \in V\}$\tcp*{new set of vertices}
$A \gets \{(\psi(u), \psi(v), i, j) : (u,v,i,j) \in A,\ \psi(u) \neq \psi(v)\}$\tcp*{relabel the graph and remove self-loops}
\Return $(G = (V, A), \vS)$
}
\end{algorithm}

\begin{definition}
A bin type $t_1$ of capacity $W_{t_1}$ dominates a bin type $t_2$ of
capacity $W_{t_2}$, $(t_1, W_{t_1})~\prec~(t_2, W_{t_2})$ for short, 
if $W_{t_1}~=~W_{t_2}$ and $t_1~<~t_2$, or $W_{t_1}~\neq~W_{t_2}$ and
$W_{t_1}~\leq~W_{t_2}$.
\end{definition}

\begin{algorithm}[!h]
\caption{Connect internal nodes to the targets}\label{alg:contgts}
\SetEndCharOfAlgoLine{}
\SetKwInOut{Input}{input}
\SetKwInOut{Output}{output}

\SetKwFunction{Build}{build}
\SetKwFunction{BuildGraph}{buildGraph}
\SetKwFunction{ConnectTargets}{connectTargets}
\SetKwFunction{Sort}{sorted}
\SetKwFunction{Key}{key}
\SetKwFunction{Len}{len}
\SetKwFunction{Reverse}{reverse}
\SetKwFunction{Lift}{lift}
\SetKwFunction{HPos}{highestPosition}
\SetKwData{varitems}{items}
\SetKwData{vardp}{dp}
\SetKwData{varx}{x}
\SetKwData{vari}{i}
\SetKwData{varmx}{mx}
\SetKwData{varmi}{mi}
\SetKwData{NIL}{NIL}

\SetInd{0.5em}{1.5em}

\SetKwBlock{Function}{}{}

\Input{$V$ - set of vertices; $A$ - set of arcs; $\vS$ - source; $q$ - number of bin types; $W$ - capacity vectors}
\Output{$V$ - set of vertices; $A$ - set of arcs; $\vS$ - source; $\vTs$ - targets}

\textbf{function} $\ConnectTargets(V, A, \vS, q, W)$:
\Function{
$\vTs \gets \langle \vTs_t \rangle_{t=1,\ldots,q}$\\
\ForEach(\tcp*[f]{for each internal node}){$v \in V \setminus \{\vS\}$}{
	$\tau \gets \{t : t =1, \ldots, q, v \leq W_t\}$\tcp*{valid bin types for vertex $v$}
	\ForEach(){$t \in \{1,\ldots,q\} $}{	
		\If(){$t \in \tau$}{
			$\tau \gets \tau \setminus \{t' \in \tau : (t, W_t) \prec (t', W_{t'})  \}$\tcp*{transitive reduction (i.e., remove dominated bin types)} 
		}
	}
	$A \gets A \cup \{(v, \vTs_t, 0, 0): t \in \tau \}$\tcp*{connect $v$ to non-dominated targets}	
}
\ForEach(\tcp*[f]{for each bin type}){$t \in \{1,\ldots,q\} $}{	
	$\tau \gets \{t' = 1,\ldots,q : (t, W_t) \prec (t', W_{t'}) \}$\tcp*{dominated bin types}
	\ForEach(){$t' \in \{1,\ldots,q\} $}{	
		\If(){$t' \in \tau$}{
			$\tau \gets \tau \setminus \{t'' \in \tau : (t', W_{t'}) \prec (t'', W_{t''})  \}$\tcp*{transitive reduction}
		}
	}
	$A \gets A \cup \{(\vTs_t, \vTs_{t'}, 0, 0): t' \in \tau \}$\tcp*{connect $\vTs_t$ the targets it directly dominates}
}
$V \gets V \cup \{\vTs_t : t=1,\ldots,q \}$\\
$A \gets A \cup \{(\vTs_t, \vS, 0, 0) : t=1,\ldots,q \}$\tcp*{add the feedback arcs}
\Return $(G = (V, A), \vS, \vTs)$
}
\end{algorithm}

\FloatBarrier
\section{Conclusions}
\label{sec:conclusions}

In this paper, we presented a graph construction and compression algorithm for multiple-choice
vector bin packing problems. This algorithm was introduced in VPSolver 3.0.0 (\url{https://github.com/fdabrandao/vpsolver}) as the standard
graph construction method. 
When there is only one bin type and each item has a single incarnation, the new algorithm
produces exactly the same graph as the one produced by the original algorithm.
When multiple bin types exist, the algorithm produces on a single run
a graph containing all the valid packing patterns for all bin types.

\FloatBarrier

\bibliographystyle{apalike} 
\bibliography{arcflow,paper}
 
\end{document}